\documentclass[11pt]{article}
\usepackage{latexsym}
\usepackage{amsfonts}
\usepackage{amssymb}

\usepackage[active]{srcltx} 
\usepackage{color}
\usepackage{hyperref}
\usepackage{refcheck}

\usepackage{amsmath}

\usepackage[active]{srcltx} 

\newcommand{\ba}{\begin{array}}\newcommand{\ea}{\end{array}}

\newcommand{\ns}{\rm}

\newcommand{\nse}{\kern-3pt\ns$=$}\newcommand{\qd}{\hfill$\Box$\medbreak}

\newcommand{\ext}{\raise1pt\hbox{$\ts\bigwedge$}}

\newcommand{\ts}{\textstyle}
\newcommand{\rf}[1]{(\ref{#1})}
\newcommand{\chii}{\raise2pt\hbox{$\chi$}}

\newcommand{\Fg}{\mbox{${\cal F}\kern-2pt_g$}}

\newcommand{\Mg}{\mbox{${\cal M}\kern-2pt_g$}}
\newcommand{\Ng}{\mbox{${\cal N}\kern-2pt_g$}}
\newcommand{\V}{V\kern-1pt}

\newcommand{\Gg}{\mbox{${\cal G}\kern-2pt_g$}}
\newcommand{\cir}{\raise1.6pt\hbox{\footnotesize$\circ$}}

\newcommand{\Res}[2]{\hbox{\ns Res}\kern-16pt\lower5pt\hbox{\footnotesize$_{#1}$}\kern2pt\left[#2\right]}
\newcommand{\qk}{quaternion-K\"ahler\kern2pt}\renewcommand{\,}{\kern1pt}

\newcommand{\End}{{\rm End}}

\newcommand{\dirac}{/\kern-5pt\partial}

\newcommand{\lra}{\longrightarrow}
\renewcommand{\ts}{\textstyle}
0
\newtheorem{theo}{Theorem}[section]
\newtheorem{defi}{Definition}[section]\newtheorem{lemma}{Lemma}[section]

\newtheorem{corol}{Corollary}[section]
\def\frac#1#2{{#1\over#2}}
\def\be#1\ee{\begin{equation}#1\end{equation}}

\allowdisplaybreaks

\begin{document}

\title{Twisted partially pure spinors}

\author{
Rafael Herrera\footnote{Centro de
Investigaci\'on en Matem\'aticas, A. P. 402,
Guanajuato, Gto., C.P. 36000, M\'exico. E-mail: rherrera@cimat.mx}
\footnote{Partially supported by 
grants of CONACyT, LAISLA (CONACyT-CNRS), and the IMU Berlin Einstein Foundation
Program} \,\, and
Ivan Tellez\footnote{Centro de
Investigaci\'on en Matem\'aticas, A. P. 402,
Guanajuato, Gto., C.P. 36000, M\'exico. E-mail: tellezito@cimat.mx}
\footnote{Partially supported by a
CONACYT scholarship}
 }

\date{}

\maketitle

\vspace{-20pt}

{
\abstract{

Motivated by the relationship between orthogonal complex structures and spure spinors,
we define twisted partially pure spinors
in order to characterize spinorially subspaces of Euclidean space endowed with a complex structure. 

}
}

\section{Introduction}

In this paper, we characterize subspaces of Euclidean space $\mathbb{R}^{n}$ endowed with an orthogonal
complex structure by
means of twisted spinors, which
is a generalization of the relation between classical pure spinors and orthogonal complex
structures on Euclidean space $\mathbb{R}^{2m}$.
Recall that a classical pure spinor $\phi\in\Delta_{2m}$ is a spinor such
that the (isotropic) subspace of complexified vectors $X-iY\in\mathbb{R}^{2m}\otimes \mathbb{C}$,
$X,Y\in\mathbb{R}^{2m}$, which
annihilate $\phi$ under Clifford multiplication
\[(X-iY)\cdot \phi =0\]
is of maximal dimension, where $m\in\mathbb{N}$ and $\Delta_{2m}$ is
the standard complex representation of the Spin group $Spin(2m)$ (cf. \cite{Lawson}). This means that
for every $X\in\mathbb{R}^{2m}$ there exists a $Y\in\mathbb{R}^{2m}$ satisfying 
\[X\cdot \phi = i Y\cdot \phi. \]
By setting $Y=J(X)$, one can see that a pure spinor determines a complex
structure on $\mathbb{R}^{2m}$.
Geometrically, the two subspaces $TM\cdot\phi$ and $i\,TM\cdot\phi$ of 
$\Delta_{2m}$ coincide, which means $TM\cdot \phi$ is a complex subspace of $\Delta_{2m}$, and
the effect of multiplication by the number $i=\sqrt{-1}$ is transferred to the tangent space $TM$ in
the form of $J$. 

The authors of \cite{Charlton,Trautman} investigated (the classification of) non-pure classical spinors 
by means of their isotropic subspaces. In \cite{Trautman}, the authors noted that there may be many spinors
(in different orbits under the action of the Spin group) admitting isotropic subspaces of the same dimension,
and that there is a gap in the possible dimensions of such isotropic subspaces.
In our Euclidean/Riemannian
context, such isotropic subspaces correspond to subspaces of
Euclidean space endowed with orthogonal complex structures. 
In this paper, we define twisted partially pure spinors (cf. Definition
\ref{def:twisted-partially-pure-spinor}) in order to  establish a one-to-one
correspondence between subspaces of Euclidean space (of a fixed codimension) endowed with orthogonal complex
structures (and oriented orthogonal complements), and orbits of such spinors under a particular subgroup
of the twisted spin group (cf. Theorem \ref{theo:characterization}). By using spinorial twists we avoid having
different orbits under the full twisted spin group and also the aforementioned gap in the dimensions.

The need to establish such a correspondence arises from our interest in developing a spinorial 
setup to study the geometry of manifolds admitting (almost) CR structures (of arbitrary codimension) and
elliptic structures. Since such manifolds are not necessarily Spin nor Spin$^c$, we are led to consider
spinorially twisted spin groups, representations, structures, etc. 
Geometric and topological considerations regarding such manifolds will be presented
in \cite{Herrera-Nakad}. 

The paper is organized as follows. In Section \ref{sec:preliminaries} we recall basic material on Clifford
algebras, spin groups and representations; we define the twisted spin groups and representations that will be
used, and the space of anti-symmetric 2-forms and endomorphims associated to twisted spinors; we also present
some results on subgroups and branching of representations.
In Section \ref{sec:twisted-partially-pure-spinors}, we define partially pure spinors, deduce their
basic properties and prove the main theorem, Theorem \ref{theo:characterization}, which establishes the
aforementioned one-to-one correspondence.

{\bf Acknowledgments}. 
The first named author would like to thank Helga Baum for her hospitality and support, as well as the
following
institutions: Humboldt University, the 
International Centre for Theoretical Physics and the Institut des Hautes \'Etudes Scientifiques.

\section{Preliminaries}\label{sec:preliminaries}

In this section, we briefly recall basic facts about  Clifford algebras, the Spin
group and the standard Spin representation \cite{Friedrich}. We also define the twisted spin groups and
representations, and the antisymmetric 2-forms and endomorphisms associated to a twisted spinor, and
describe various inclusions of groups into (twisted) spin groups.

\subsection{Clifford algebras}

Let $Cl_n$ denote the Clifford algebra generated by the orthonormal vectors
$e_1, e_2, \ldots, e_n\in \mathbb{R}^n$ 
subject to the relations
\begin{eqnarray*}
e_j e_k + e_k e_j &=& -2\left< e_j,e_k\right>,
\end{eqnarray*}
where $\big< , \big>$ denotes the standard inner product in $\mathbb{R}^n$.
Let
\[\mathbb{C}l_n=Cl_n\otimes_{\mathbb{R}}\mathbb{C}\]
denote the complexification of $Cl_n$. The Clifford algebras are isomorphic to matrix algebras. In particular,
\[\mathbb{C}l_n\cong \left\{
                     \begin{array}{ll}
                     \End(\mathbb{C}^{2^k}), & \mbox{if $n=2k$,}\\
                     \End(\mathbb{C}^{2^k})\oplus\End(\mathbb{C}^{2^k}), & \mbox{if $n=2k+1,$}
                     \end{array}
\right.
\]
where
\[\Delta_n:=\mathbb{C}^{2^k}=\underbrace{\mathbb{C}^2\otimes \ldots \otimes \mathbb{C}^2}_{k\,\,\,\,
{\rm times}}\]
is the tensor product of $k=[{n\over 2}]$ copies of $\mathbb{C}^2$.
The map
\[\kappa:\mathbb{C}l_n \lra \End(\mathbb{C}^{2^k})\]
is defined to be either the above mentioned isomorphism if $n$ is even, or the isomorphism followed by
the projection onto the first summand if $n$ is odd.
In order to make $\kappa$ explicit, consider the following matrices
\[Id = \left(\begin{array}{ll}
1 & 0\\
0 & 1
      \end{array}\right),\quad
g_1 = \left(\begin{array}{ll}
i & 0\\
0 & -i
      \end{array}\right),\quad
g_2 = \left(\begin{array}{ll}
0 & i\\
i & 0
      \end{array}\right),\quad
T = \left(\begin{array}{ll}
0 & -i\\
i & 0
      \end{array}\right).
\]
In terms of the generators $e_1, \ldots, e_n$, $\kappa$ can be
described explicitly as follows,
\begin{eqnarray}
e_1&\mapsto& Id\otimes Id\otimes \ldots\otimes Id\otimes Id\otimes g_1,\nonumber\\
e_2&\mapsto& Id\otimes Id\otimes \ldots\otimes Id\otimes Id\otimes g_2,\nonumber\\
e_3&\mapsto& Id\otimes Id\otimes \ldots\otimes Id\otimes g_1\otimes T,\nonumber\\
e_4&\mapsto& Id\otimes Id\otimes \ldots\otimes Id\otimes g_2\otimes T,\nonumber\\
\vdots && \dots\nonumber\\
e_{2k-1}&\mapsto& g_1\otimes T\otimes \ldots\otimes T\otimes T\otimes T,\nonumber\\
e_{2k}&\mapsto& g_2\otimes T\otimes\ldots\otimes T\otimes T\otimes T,\nonumber
\end{eqnarray}
and, if $n=2k+1$, 
\[ e_{2k+1}\mapsto i\,\, T\otimes T\otimes\ldots\otimes T\otimes T\otimes T.\]
The vectors 
\[u_{+1}={1\over \sqrt{2}}(1,-i)\quad\quad\mbox{and}\quad\quad u_{-1}={1\over \sqrt{2}}(1,i),\]
form a unitary basis of $\mathbb{C}^2$ with respect to the standard Hermitian product.
Thus, 
\[\{u_{\varepsilon_1,\ldots,\varepsilon_k}=u_{\varepsilon_1}\otimes\ldots\otimes
u_{\varepsilon_k}\,\,|\,\, \varepsilon_j=\pm 1,
j=1,\ldots,k\},\]
is a unitary basis of $\Delta_n=\mathbb{C}^{2^k}$
with respect to the naturally induced Hermitian product.
We will denote inner and Hermitian products  by the same symbol $\left<\cdot,\cdot\right>$ trusting that the
context will make clear
which product is being used.

Clifford multiplication is defined by
\begin{eqnarray*}
\mu_n:\mathbb{R}^n\otimes \Delta_n &\lra&\Delta_n\\ 
x \otimes \psi &\mapsto& \mu_n(x\otimes \psi)=x\cdot\psi :=\kappa(x)(\psi) .
\end{eqnarray*}

A quaternionic structure $\alpha$ on $\mathbb{C}^2$ is given by
\[\alpha\left(\begin{array}{c}
z_1\\
z_2
              \end{array}
\right) = \left(\begin{array}{c}
-\overline{z}_2\\
\overline{z}_1
              \end{array}\right),\]
and a real structure $\beta$ on $\mathbb{C}^2$ is given by
\[\beta\left(\begin{array}{c}
z_1\\
z_2
              \end{array}
\right) = \left(\begin{array}{c}
\overline{z}_1\\
\overline{z}_2
              \end{array}\right).\]
Real and quaternionic structures $\gamma_n$  on $\Delta_n=(\mathbb{C}^2)^{\otimes
[n/2]}$ are built as follows
\[
\begin{array}{cclll}
 \gamma_n &=& (\alpha\otimes\beta)^{\otimes 2k} &\mbox{if $n=8k,8k+1$}& \mbox{(real),} \\
 \gamma_n &=& \alpha\otimes(\beta\otimes\alpha)^{\otimes 2k} &\mbox{if $n=8k+2,8k+3$}&
\mbox{(quaternionic),} \\
 \gamma_n &=& (\alpha\otimes\beta)^{\otimes 2k+1} &\mbox{if $n=8k+4,8k+5$}&\mbox{(quaternionic),} \\
 \gamma_n &=& \alpha\otimes(\beta\otimes\alpha)^{\otimes 2k+1} &\mbox{if $n=8k+6,8k+7$}&\mbox{(real).}
\end{array}
\]

\subsection{The Spin group and representation}
The Spin group $Spin(n)\subset Cl_n$ is the subset 
\[Spin(n) =\{x_1x_2\cdots x_{2l-1}x_{2l}\,\,|\,\,x_j\in\mathbb{R}^n, \,\,
|x_j|=1,\,\,l\in\mathbb{N}\},\]
endowed with the product of the Clifford algebra.
It is a Lie group and its Lie algebra is
\[\mathfrak{spin}(n)=\mbox{span}\{e_ie_j\,\,|\,\,1\leq i< j \leq n\}.\]
Recall that the Spin group $Spin(n)$ is the universal double cover of $SO(n)$, $n\ge 3$. For $n=2$
we consider $Spin(2)$ to be the connected double cover of $SO(2)$.
The covering map will be denoted by 
\[\lambda_n:Spin(n)\rightarrow SO(n).\]
Its differential is given
by $\lambda_{n_*}(e_ie_j) = 2E_{ij}$, where $E_{ij}=e_i^*\otimes e_j - e_j^*\otimes e_i$ is the
standard basis of the skew-symmetric matrices, and $e^*$ denotes the metric dual of the vector $e$.
Furthermore, we will abuse the notation and also denote by $\lambda_n$ the induced representation on
$\ext^*\mathbb{R}^n$.

The restriction of $\kappa$ to $Spin(n)$ defines the Lie group representation
\[
Spin(n)\lra GL(\Delta_n),\]
which is, in fact, special unitary. 
We have the corresponding Lie algebra representation
\[
\mathfrak{spin}(n)\lra \mathfrak{gl}(\Delta_n).\] 

{\bf Remark}.
For the sake of notation we will set
\[SO(0)=\{1\},\quad\quad SO(1)=\{1\},\]
\[Spin(0)=\{\pm1\},\quad\quad Spin(1)=\{\pm1\},\]
and
\[\Delta_0 = \Delta_1 =\mathbb{C}\]
a trivial $1$-dimensional representation.

Clifford multiplication $\mu_n$ has the following properties:
\begin{itemize}
 \item It is skew-symmetric with respect to the Hermitian product
\begin{equation}
\left<x\cdot\psi_1 , \psi_2\right> 
=-\left<\psi_1 , x\cdot \psi_2\right>. \label{clifford-skew-symmetric}
\end{equation} 
 \item $\mu_n$ is an equivariant map of $Spin(n)$ representations.
 \item $\mu_n$ can be extended to an equivariant map 
\begin{eqnarray*}
\mu_n:\ext^*(\mathbb{R}^n)\otimes \Delta_n &\lra&\Delta_n\\ 
\omega \otimes \psi &\mapsto& \omega\cdot\psi,
\end{eqnarray*}
of $Spin(n)$ representations. 
\end{itemize}

At this point we will make the following convention. Consider the involution 
\begin{eqnarray*}
F_{2m}: \Delta_{2m}&\lra& \Delta_{2m}\\
\phi &\mapsto& (-i)^m e_1e_2\cdots e_{2m} \cdot \phi,
\end{eqnarray*}
and let
\[\Delta_{2m}^\pm = \{\phi \,\,|\,\, F_{2m}(\phi)=\pm \phi \}.\]
This definition of positive and negative Weyl spinors differs from the one in \cite{Friedrich} by a factor
$(-1)^m$. Nevertheless, we have chosen this convention so that the spinor $u_{1,\ldots,1}$ is
always
positive and corresponds to the standard (positive) complex structure on $\mathbb{R}^{2m}$.

\subsection{Spinorially twisted Spin groups}

Consider the following groups:
\begin{itemize}
 \item By using the unit complex numbers $U(1)$,
the Spin group can be twisted \cite{Friedrich}
\[Spin^c(n) =  (Spin(n) \times U(1))/\{\pm (1,1)\} =
Spin(n) \times_{\mathbb{Z}_2} U(1),\]
with Lie algebra
\[\mathfrak{spin}^c(n)=\mathfrak{spin}(n)\oplus i\mathbb{R}.\]

\item In \cite{Espinosa-Herrera} we have considered
 the twisted Spin group $Spin^r(n)$, $r\in\mathbb{N}$, defined as follows 
\[
Spin^{r}(n) =  (Spin(n) \times Spin(r))/\{\pm (1,1)\} =
Spin(n) \times_{\mathbb{Z}_2} Spin(r).
\]
The Lie algebra of $Spin^r(n)$  is
\[\mathfrak{spin}^r(n) = \mathfrak{spin}(n) \oplus \mathfrak{spin}(r).\]

\item Here, we will also consider the following group
\begin{eqnarray*}
Spin^{c,r}(n) &=&  (Spin(n) \times Spin^c(r))/\{\pm (1,1)\} \\
&=& Spin(n) \times_{\mathbb{Z}_2} Spin^c(r),
\end{eqnarray*}
where $r\in\mathbb{N}$, whose
Lie algebra is
\[\mathfrak{spin}^c(n)=\mathfrak{spin}(n)\oplus \mathfrak{spin}(r)\oplus i\mathbb{R}.\]
It fits into the exact sequence
\[1\lra \mathbb{Z}_2\lra Spin^{c,r}(n)\xrightarrow{\lambda_n\times\lambda_r\times \lambda_2} SO(n)\times SO(r)
\times U(1)\lra 1,\]
where
\begin{eqnarray*}
(\lambda_n\times\lambda_r\times \lambda_2)([g,[h,z]])  &=& (\lambda_n(g),\lambda_r(h),z^2).
\end{eqnarray*}
\end{itemize}

{\bf Remark}. For $r=0,1$,  $Spin^{c,r}(n)=Spin^c(n)$.

\subsection{Twisted spin representations}

Consider the following twisted representations:
\begin{itemize}
 \item The Spin representation $\Delta_n$ extends to a representation of $Spin^c(n)$ by letting
\begin{eqnarray*}
Spin^c(n)&\longrightarrow& GL(\Delta_n)\\
\,[g,z]  &\mapsto& z\kappa_n(g)=:zg.
\end{eqnarray*}

\item The twisted $Spin^{c,r}(n)$ representation
\begin{eqnarray*}
Spin^{c,r}(n)&\longrightarrow& GL(\Delta_r\otimes
\Delta_n)\\
\,[g,[h,z]]  &\mapsto& z \, \kappa_r(h)\otimes \kappa_n(g)=:zh\otimes g.
\end{eqnarray*}
which is also unitary with respect to the natural Hermitian metric.
\item For $r=0,1$, the twisted spin representation is simply the $Spin^c(n)$ representation $\Delta_n$.
\end{itemize}

We will also need the map
\begin{eqnarray*}
 \mu_r\otimes\mu_n:\left(\ext^*\mathbb{R}^r\otimes_\mathbb{R} \ext^*\mathbb{R}^n\right)
 \otimes_\mathbb{R} (\Delta_r\otimes \Delta_n) &\longrightarrow& \Delta_r\otimes\Delta_n\\
(w_1 \otimes w_2)\otimes (\psi\otimes \varphi) &\mapsto& 
(w_1\otimes w_2)\cdot (\psi\otimes \varphi)
= (w_1\cdot\psi) \otimes (w_2\cdot \varphi).
\end{eqnarray*}
As in the untwisted case, 
$\mu_r\otimes\mu_n$ is an equivariant homomorphism of $Spin^{c,r}(n)$ representations.

\subsection{Skew-symmetric 2-forms and endomorphisms associated to twisted spinors}

We will often write $f_{kl}$ for the Clifford product $f_kf_l$.

\begin{defi}
{\rm \cite{Espinosa-Herrera}}
Let $r\geq 2$, $\phi\in\Delta_r\otimes\Delta_n$, $X,Y\in\mathbb{R}^n$, $(f_1\ldots,f_r)$
an orthonormal basis of $\mathbb{R}^r$ and $1\leq k,l\leq r$.
\begin{itemize}
\item 
Define the real $2$-forms associated to the spinor $\phi$ by
\[\eta_{kl}^{\phi} (X,Y) = {\rm Re}\left< X\wedge Y\cdot f_kf_l\cdot \phi,\phi\right>.\]

\item Define the antisymmetric endomorphisms
$\hat\eta_{kl}^\phi\in\End^-(\mathbb{R}^n)$ by
\[X\mapsto \hat\eta_{kl}^\phi(X):=(X\lrcorner \,\eta_{kl}^{\phi})^\sharp,\]
where $X\in\mathbb{R}^n$, $\lrcorner$ denotes contraction and $^\sharp$ denotes metric dualization from
$1$-fomrs to vectors.
\end{itemize}
\end{defi}

\begin{lemma}
Let $r\geq 2$, $\phi\in\Delta_r\otimes\Delta_n$, $X,Y\in\mathbb{R}^n$, $(f_1\ldots,f_r)$
an orthonormal basis of $\mathbb{R}^r$ and $1\leq k,l\leq r$. Then
\begin{eqnarray}
{\rm Re}\left<  f_kf_l\cdot
\phi,\phi\right>&=&0,\label{vanishing2}\nonumber\\
{\rm Re}\left< X\wedge Y\cdot
\phi,\phi\right>&=&0,\label{vanishing3}\\
{\rm Im}\left< X\wedge Y\cdot f_kf_l\cdot
\phi,\phi\right>&=&0,\label{vanishing4}\\
{\rm Re} \left< X\cdot \phi,Y\cdot\phi \right> 
 &=&   \left<X,Y\right>|\phi|^2, \label{real-part}
\end{eqnarray}
\end{lemma}
{\em Proof}.
By using \rf{clifford-skew-symmetric} twice
\begin{eqnarray*}
 \left<  f_kf_l\cdot\phi,\phi\right>
&=& -\overline{\left<   f_kf_l\phi,\phi\right>}.
\end{eqnarray*}

For identity \rf{vanishing3}, recall that for $X, Y\in \mathbb{R}^n$ 
\[X\wedge Y =  X\cdot Y + \left<X,Y\right>.\]
Thus
\begin{eqnarray*}
 \left< X\wedge Y\cdot \phi,\phi\right>
 &=& -\overline{\left< X\wedge Y \cdot\phi,\phi\right>}.
\end{eqnarray*}

Identities \rf{vanishing4} and \rf{real-part} follow similarly.
\qd

{\bf Remarks}. 
\begin{itemize}
 \item For $k\not= l$,
\[\eta_{kl}^\phi =  (\delta_{kl}-1)\eta_{lk}^\phi.\]
 \item By \rf{vanishing4}, if $k\not= l$,
\[\eta_{kl}^{\phi} (X,Y) =\left< X\wedge Y\cdot f_kf_l\cdot
\phi,\phi\right>.\]
\end{itemize}

\begin{lemma}{\rm \cite{Espinosa-Herrera}}
Any spinor $\phi\in\Delta_r\otimes\Delta_n$, $r\geq 2$, defines two maps (extended by linearity)
\begin{eqnarray*}
\ext^2 \mathbb{R}^r&\lra& \ext^2 \mathbb{R}^n\\
f_{kl} &\mapsto& \eta_{kl}^{\phi}
\end{eqnarray*}
and
\begin{eqnarray*}
\ext^2 \mathbb{R}^r&\lra& \End(\mathbb{R}^n)\\
f_{kl} &\mapsto& \hat\eta_{kl}^{\phi},
\end{eqnarray*}
\end{lemma}
\qd

\subsection{Subgroups, isomorphisms and decompositions}

In this section we will describe various inclusions of groups into (twisted) spin groups.

\begin{lemma}\label{lemma-subgroup1}
There exists a monomorphism $h:Spin(2m)\times_{\mathbb{Z}_2} Spin(r) \lra Spin(2m+r)$
such that the following diagram commutes
\[
 \begin{array}{ccc}
Spin(2m)\times_{\mathbb{Z}_2} Spin(r) & \xrightarrow{h} & Spin(2m+r)\\
\downarrow &  & \downarrow\\
SO(2m)\times SO(r) & \hookrightarrow & SO(2m+r)
 \end{array}
\]
\end{lemma}
{\em Proof}.
Consider the decomposition
\[\mathbb{R}^{2m+r}=\mathbb{R}^{2m}\oplus\mathbb{R}^{r},\]
and let
\begin{eqnarray*}
 Spin(2m) &=& \left\{\prod_{i=1}^{2s} x_i \in Cl_{2m+r} \,\, | \,\,
 x_i\in\mathbb{R}^{2m}, |x_i|=1, s\in\mathbb{N} \right\}\quad\subset\quad Spin(2m+r),\\
 Spin(r) &=& \left\{\prod_{j=1}^{2t} y_j \in Cl_{2m+r} \,\, | \,\,
 y_j\in\mathbb{R}^{r}, |y_j|=1, t\in\mathbb{N} \right\}\quad\subset\quad Spin(2m+r).
\end{eqnarray*}
It is easy to prove that
\[Spin(2m)\cap Spin(r) = \{1,-1\}.\]

Define the homomorphism
\begin{eqnarray*}
 h:Spin(2m)\times_{\mathbb{Z}_2} Spin(r) &\lra& Spin(2m+r)\\
{[g,g']} &\mapsto & gg'.
\end{eqnarray*}
If $[g,g']\in Spin(2m)\times_{\mathbb{Z}_2} Spin(r)$ is such that
\[gg'=1\in Spin(2m+r),\]
then
\[g'=g^{-1} \in Spin(2m)\subset Spin(2m+r),\]
so that 
\[g,g'\in Spin(2m)\cap Spin(r)=\{1,-1\}.\]
Hence $[g,g']=[1,1]$
and $h$ is injective.
\qd

\begin{lemma}\label{lemma:subgroup2}
Let $r\in \mathbb{N}$.
 There exists an monomorphism $h:U(m)\times SO(r) \hookrightarrow
Spin^{c,r}(2m+r)$
such that the following diagram commutes
\[
 \begin{array}{ccc}
 &  & Spin^{c,r}(2m+r)\\
 & \nearrow & \downarrow\\
U(m)\times SO(r) & \lra & SO(2m+r)\times SO(r)\times U(1)
 \end{array}
\]
\end{lemma}
{\em Proof}.
Suppose we have an orthogonal complex structure on $\mathbb{R}^{2m}\subset\mathbb{R}^{2m+r}$
\[J:\mathbb{R}^{2m} \lra \mathbb{R}^{2m}, \quad J^2=-{\rm Id}_{2m}, \quad
\langle\cdot,\cdot\rangle=\langle J\cdot,J\cdot\rangle.\]
The subgroup of $SO(2m+r)$ that respects both the orthogonal decomposition 
$\mathbb{R}^{2m+r}=\mathbb{R}^{2m}\oplus\mathbb{R}^{r}$
and $J$ is
\[U(m)\times SO(r)\quad\subset\quad SO(2m)\times SO(r)\quad\subset\quad SO(2m+r).\]
There exists a lift \cite{Friedrich}
\[\begin{array}{ccc}
 &  & Spin^c(2m)\\
 & \nearrow & \downarrow\\
U(m) & \rightarrow & SO(2m)\times U(1)\\
&&\\
A & \mapsto &(A_{\mathbb{R}},\det_{\mathbb{C}}(A))
  \end{array}
\]
and we can consider the diagram \cite{Espinosa-Herrera}
\[
 \begin{array}{ccc}
 &  & Spin(r)\times_{\mathbb{Z}_2} Spin(r)\\
 & \nearrow & \downarrow\\
SO(r) & \xrightarrow{\rm diagonal} & SO(r)\times SO(r)
 \end{array}.
\]
We can put them together as follows
{\footnotesize
\[
 \begin{array}{ccccccc}
 & & Spin^c(2m)\times_{\mathbb{Z}_2} Spin^r(r) \cong  Spin^r(2m)\times_{\mathbb{Z}_2}
Spin^c(r)&\hookrightarrow& Spin(2m+r)\times_{\mathbb{Z}_2} Spin^c(r)\\
 &\nearrow & \downarrow & &\\
U(m)\times SO(r) & \hookrightarrow & SO(2m)\times U(1)\times SO(r)\times SO(r)
         &&
 \end{array}
\]
}
where the last inclusion is due to Lemma \ref{lemma-subgroup1}.
It is easy to prove that 
the lift monomorphism $U(m)\times SO(r)\lra Spin^c(2m)\times_{\mathbb{Z}_2} Spin^r(r)$ 
exists and 
there is a natural isomorphism
\[Spin^c(2m)\times_{\mathbb{Z}_2} Spin^r(r) \cong  Spin^r(2m)\times_{\mathbb{Z}_2} Spin^c(r).\]
\qd

\begin{lemma}\label{factorization}
Let $r\in\mathbb{N}$.
 The standard representation $\Delta_{2m+r}$ of $Spin(2m+r)$ decomposes as follows
\[\Delta_{2m+r} = \Delta_r\otimes\Delta_{2m}^+ \,\,\oplus\,\,\Delta_r\otimes\Delta_{2m}^-,\]
with respect to the subgroup
$Spin(2m)\times_{\mathbb{Z}_2} Spin(r) \subset Spin(2m+r)$. 
\end{lemma}
{\em Proof}. Consider the restriction of the standard representation of $Spin(2m+r)$ to
\[Spin(2m)\times_{\mathbb{Z}_2} Spin(r)\subset Spin(2m+r) \lra Gl(\Delta_{2m+r}).\]
By using the explicit description of a unitary basis of $\Delta_{2m+r}$, we see that
the elements of $Spin(2m)$ act on the last $m$ factors of
\[\Delta_{2m+r} = \underbrace{\mathbb{C}^2 \otimes \cdots\otimes  \mathbb{C}^2}_{[r/2] \: {\rm times}} \otimes
\underbrace{\mathbb{C}^2 \otimes\cdots\otimes\mathbb{C}^2 }_{m \:
{\rm times}}, \]
as they do on $\Delta_{2m}=\Delta_{2m}^+\oplus\Delta_{2m}^-$. The elements of $Spin(r)$ 
act as usual on the first $[r/2]$ factors of $\Delta_r$,
act trivially on $\Delta_{2m}^+$, 
and act by multiplication by $(-1)$ on the factor $\Delta_{2m}^-$.
\qd

\section{Twisted partially pure spinors}\label{sec:twisted-partially-pure-spinors}

In order to simplify the statements, we will
consider the twisted spin representation
\[
 \Sigma_r\otimes\Delta_n \subseteq \Delta_r\otimes\Delta_n.
\]
where 
\[
\Sigma_r=
\left\{
\begin{array}{ll}
\Delta_r & \mbox{if $r$ is odd,}\\
\Delta_r^+ & \mbox{if $r$ is even,} 
\end{array}
\right. 
\]
$n,r\in\mathbb{N}$.

\begin{defi}\label{def:twisted-partially-pure-spinor}
Let $(f_1,\ldots,f_r)$ be an  orthonormal frame of $\mathbb{R}^r$.
A unit-length spinor $\phi\in\Sigma_r\otimes\Delta_n$, $r<n$, is called a {\em twisted partially pure
spinor} if 
\begin{itemize}
 \item there exists a $(n-r)$-dimensional subspace $V^\phi\subset\mathbb{R}^n$ such that
for every $X\in V^\phi$, there exists a $Y\in V^\phi$ such that
\[X\cdot \phi = i\,\,Y\cdot \phi. \]

 \item it satisfies the equations
\begin{eqnarray*}
(\eta_{kl}^\phi + f_kf_l)\cdot \phi&=&0,\\
 \left<f_kf_l\cdot \phi,\phi\right>&=&0,
\end{eqnarray*}
for all $1\leq k<l\leq r$.
\item If $r=4$, it also satisfies the condition
\[\left<f_1f_2f_3f_4\cdot \phi,\phi\right>=0.\]

\end{itemize}
\end{defi}

{\bf Remarks}.
\begin{enumerate}

 \item The requirement $|\phi|=1$ is made in order to avoid renormalizations later on.

 \item The extra condition for the case $r=4$ is fulfilled for all other ranks.

 \item From now on we will drop the adjective twisted since it will be clear from the context.
\end{enumerate}

\subsection{Example of partially pure spinor}

\begin{lemma}\label{lemma:existence}
 Given $r,m\in\mathbb{N}$, there exists a partially pure spinor in
$\Sigma_r\otimes\Delta_{2m+r}$.
\end{lemma}
{\em Proof}.
Let $(e_1,\ldots,e_{2m},e_{2m+1},\ldots,e_{2m+r})$ and 
$(f_1,\ldots,f_r)$ be orthonormal frames 
of $\mathbb{R}^{2m+r}$ and $\mathbb{R}^r$ respectively.
Consider the decomposition of Lemma \ref{factorization}
\[\Delta_{2m+r}=\Delta_r\otimes\Delta_{2m}^+\,\,\oplus\,\,\Delta_r\otimes\Delta_{2m}^-,\]
corresponding to the decomposition
\[\mathbb{R}^{2m+r}={\rm span}\{e_1,\ldots,e_{2m}\}
\oplus {\rm span}\{e_{2m+1},\ldots,e_{2m+r}\}.\]
Let
\[\varphi_0= u_{1,\ldots,1}\in \Delta_{2m}^+,\]
and
\[\{v_{\varepsilon_1,\ldots,\varepsilon_{[r/2]}}\,|\, (\varepsilon_1,\ldots,\varepsilon_{[r/2]}) \in
\{\pm1\}^{[r/2]}\}\]
be the unitary basis of the twisting factor $\Delta_r=\Delta({\rm span}(f_1,\ldots,f_r))$ which contains
$\Sigma_r$.
Let us define the standard twisted partially pure spinor
$\phi_0\in\Sigma_r\otimes\Delta_r\otimes\Delta_{2m}^+$ by
\begin{equation}
\phi_0 = \left\{
\begin{array}{ll}
{1\over \sqrt{2^{[r/2]}}}\,\, \left(\sum_{I\in\{\pm1\}^{\times [r/2]}}
v_I\otimes
\gamma_{r}(u_I)\right)\otimes \varphi_0 & \mbox{if $r$ is
odd,}\\
{1\over \sqrt{2^{[r/2]-1}}}\,\, \left(\sum_{I\in\left[\{\pm1\}^{\times [r/2]}\right]_+}
v_I\otimes
\gamma_{r}(u_I)\right)\otimes \varphi_0 &\mbox{if $r$ is even}, 
\end{array} \right.  \nonumber
\end{equation}
where the elements of $\left[\{\pm1\}^{\times [r/2]}\right]_+$ 
contain an even number of $(-1)$.

Checking the conditions in the definition of partially pure spinor for $\phi_0$ is done by a (long) direct
calculation as in \cite{Espinosa-Herrera}.
For instance, taking $n=7$, $r=3$, we have 
$$\phi_0=\frac{1}{\sqrt{2}}(v_{1}\otimes\gamma_3(u_{1})\otimes u_{1}\otimes 
u_{1}+v_{-1}\otimes\gamma_3(u_{-1})\otimes u_{1}\otimes 
u_{1})$$
where $\gamma_3$ is a quaternionic structure. We check that this $\phi_0$ is a partially pure spinor. Putting 
$C=1/\sqrt{2}$ and remembering that $\gamma_3(u_{\epsilon})=-i\epsilon 
u_{-\epsilon}$, we get 
$$\phi_0=iC(v_{-1}\otimes u_{1}\otimes u_{1}\otimes 
u_{1}-v_{1}\otimes u_{-1}\otimes u_{1}\otimes 
u_{1}),$$
which has unit length.
Let $\{e_i\}$ be the standard basis of 
$\mathbb{R}^7$, so that
\begin{align*}
e_1\cdot\phi_0 & = iC(v_{-1}\otimes u_{1}\otimes u_{1}\otimes 
g_1(u_{1})-v_{1}\otimes u_{-1}\otimes u_{1}\otimes 
g_1(u_{1}))\\
& = ie_2\cdot \phi_0,
\end{align*}
and, similarly,
\begin{align*}
e_3\cdot\phi_0 
& = ie_4\cdot \phi_0.
\end{align*}
So, $\phi_0$ induces the standard complex structure on 
$V^{\phi_0}=\langle e_1,\,e_2,\,e_3,\,e_4\rangle$.
Let $\{f_i\}$ be the 
standard basis of $\mathbb{R}^3$. 
Similar calculations give
\[\eta_{kl}^{\phi_0}=e_{4+k}\wedge e_{4+l},\]
\[(\eta_{kl}^{\phi_0}+f_{kl})\cdot\phi_0=0,\] 
and
\[\left<f_{kl}\cdot\phi_0,\phi_0\right>=0.\]
\qd

\subsection{Properties of partially pure spinors}\label{sec: basic properties}

\begin{lemma}
 The definition of partially pure spinor does not depend on the choice of orthonormal basis of $\mathbb{R}^r$.
\end{lemma}
{\em Proof}.
If $r=0,1$, a partially pure spinor is a classical pure spinor for $n$ even or the straightforward
generalization of pure spinor for $n$ odd \cite[p. 336]{Lawson}.
Suppose $(f_1',\ldots, f_r')$ is another orthonormal frame of
$\mathbb{R}^r$, then
\[f_i'=\alpha_{i1}f_1+\cdots + \alpha_{ir}f_r,\]
so that the matrix $A=(\alpha_{ij})\in SO(r)$.
Let us denote
\[\eta_{kl}'^{\phi} (X,Y):={\rm Re}\left<X\wedge Y\cdot f_k'f_l'\cdot\phi,\phi\right>\]
Thus,
\begin{eqnarray*}
 \eta_{kl}'^{\phi}\cdot\phi 
   &=&  \sum_{1\leq a<b\leq n}\eta_{kl}'^{\phi}(e_a,e_b)e_ae_b
 \cdot\phi\\
   &=& \sum_{1\leq a<b\leq n}{\rm Re}\left< e_ae_b
     \cdot\left(\sum_{s=1}^r \alpha_{ks}f_s\right)\left(\sum_{t=1}^r
     \alpha_{lt}f_t\right)\cdot\phi,\phi\right>e_ae_b\cdot\phi\\
   &=& \sum_{1\leq a<b\leq n}\sum_{s=1}^r \sum_{t=1}^r
     \alpha_{ks}\alpha_{lt}{\rm Re}\left< e_ae_b
     \cdot f_sf_t\cdot\phi,\phi\right>e_ae_b\cdot\phi\\
   &=& \sum_{1\leq a<b\leq n}\sum_{s=1}^r \sum_{t=1}^r
     \alpha_{ks}\alpha_{lt}\eta_{st}^\phi(e_a,e_b) e_ae_b\cdot\phi\\
   &=& \sum_{s=1}^r \sum_{t=1}^r
     \alpha_{ks}\alpha_{lt}\eta_{st}^\phi\cdot\phi\\
   &=& -\sum_{s=1}^r \sum_{t=1}^r
     \alpha_{ks}\alpha_{lt}f_sf_t\cdot\phi\\
   &=& -\left(\sum_{s=1}^r \alpha_{ks}f_s\right)\left(\sum_{t=1}^r\alpha_{lt}f_t\right)\cdot\phi\\
   &=& -f_k'f_l'\cdot\phi.
\end{eqnarray*}

For the third part of the definition, note that
\begin{eqnarray*}
\left<f_k'f_l'\cdot\phi , \phi \right>
   &=& 
\left<
\left(\sum_{s=1}^r \alpha_{ks}f_s\right)
\left(\sum_{t=1}^r\alpha_{lt}f_t\right)\cdot\phi,\phi\right>\\
   &=& \sum_{s=1}^r \sum_{t=1}^r
     \alpha_{ks}\alpha_{lt}\left<f_sf_t\cdot\phi,\phi\right>\\
   &=& 0.
\end{eqnarray*}

For $r=4$, 
the volume form is invariant under $SO(4)$, 
$ f_1'f_2'f_3'f_4'
=f_1f_2f_3f_4$,
and
\[
\left<f_1'f_2'f_3'f_4'\cdot \phi,\phi\right>
    =  
\left<f_1f_2f_3f_4\cdot \phi,\phi\right>\\
   =
 0.
\]
\qd

\begin{lemma}
Given a partially pure spinor
$\phi\in\Sigma_r\otimes\Delta_n$, 
there exists an orthogonal complex structure on $V^\phi$ and $n-r \equiv 0$ {\rm (mod 2)} .
\end{lemma}
{\em Proof}.
By definition, for every $X\in V^\phi$, there exists $Y\in V^\phi$ such that
\[X\cdot\phi = i \, Y\cdot \phi,\]
and
\[Y\cdot\phi = i\, (-X)\cdot\phi.\] 
If we set
\[J^\phi(X) :=Y,\]
we get a linear transformation 
$J^\phi:V^\phi\rightarrow V^\phi$, such that $(J^\phi)^2=-{\rm Id}_{V^\phi}$, 
i.e. $J^\phi$ is a complex structure on the vector space $V^{\phi}$ and
$\dim_{\mathbb{R}}(V^\phi)$ is even. 
Furthermore, this  complex structure is orthogonal. Indeed, for every $X\in V^\phi$,
\begin{eqnarray*}
X\cdot JX\cdot\phi &=& -i \vert X\vert^2\phi,\\
 JX\cdot X\cdot\phi &=& i \vert JX\vert^2\phi, 
\end{eqnarray*}
and
\[  (-2 \left<X, JX\right> + i( 
\vert JX\vert^2-\vert X\vert^2))\,\,\phi=0,\]
i.e. 
\begin{eqnarray*}
\left<X, JX\right> &=&0 \\ 
\vert X\vert &=& \vert JX\vert.
\end{eqnarray*} 
\qd

\begin{lemma} Let $r\geq2$ and $\phi\in\Sigma_r\otimes\Delta_n$ be a
partially pure spinor.
The forms $\eta_{kl}^\phi$ are non-zero, $1\leq k<l\leq r$. 
\end{lemma}
{\em Proof}.
Since $(f_kf_l)^2=-1$, the equation
\begin{equation}
\eta_{kl}^\phi\cdot\phi=-f_kf_l\cdot \phi\label{eq:despeje1} 
\end{equation}
implies
\begin{equation}
\eta_{kl}^\phi\cdot f_kf_l\cdot \phi=\phi.\label{eq:despeje2} 
\end{equation}
By taking an orthonormal frame $(e_1,\ldots,e_n)$ of $\mathbb{R}^n$ we can write
\[\eta_{kl}^\phi = \sum_{1\leq i<j\leq n} \eta_{kl}^\phi(e_i,e_j)e_ie_j.\]
By \rf{eq:despeje2}, and taking hermitian product with $\phi$ 
\begin{eqnarray*}
1&=& |\phi|^2 \\
&=& \left< \eta_{kl}^\phi\cdot f_kf_l\cdot \phi,\phi\right> \\
&=& \left< \sum_{1\leq i<j\leq n} \eta_{kl}^\phi(e_i,e_j)e_ie_j\cdot f_kf_l\cdot \phi,\phi\right> \\
&=& \sum_{1\leq i<j\leq n} \eta_{kl}^\phi(e_i,e_j)\left< e_ie_j\cdot f_kf_l\cdot \phi,\phi\right> \\
&=& \sum_{1\leq i<j\leq n} \eta_{kl}^\phi(e_i,e_j)^2.
\end{eqnarray*}
\qd

\begin{lemma}\label{lemma:lie-algebra}
Let $r\geq2$.
 The image of the map associated to a partially pure spinor $\phi\in\Sigma_r\otimes\Delta_n$,
\begin{eqnarray*}
\ext^2 \mathbb{R}^r&\lra& \End(\mathbb{R}^n)\\
f_{kl} &\mapsto& \hat\eta_{kl}^{\phi},
\end{eqnarray*}
forms a Lie algebra of endomorphisms isomorphic to $\mathfrak{so}(r)$.
\end{lemma}
{\em Proof}. 
Let $(e_1,\ldots,e_n)$ be an orthonormal frame of $\mathbb{R}^n$.
First, let us consider the following calculation for $i\not=j$, $k\not =l$, $s\not=t$:
\begin{eqnarray*}
{\rm Re}\left<e_s e_t \cdot \eta_{ij}^\phi\cdot f_kf_l\cdot\phi,\phi\right> 
&=& {\rm Re}\left<e_s e_t \cdot \left(\sum_{a<b}\eta_{ij}^\phi(e_a,e_b)e_a e_b\right)\cdot
f_kf_l\cdot\phi,\phi\right>\\
&=& {\rm Re}\sum_{a<b}\eta_{ij}^\phi(e_a,e_b)\left<e_s\cdot e_t \cdot e_a\cdot e_b\cdot
f_kf_l\cdot\phi,\phi\right>\\
&=& -\sum_{b}\eta_{ij}^\phi(e_s,e_b)\eta_{kl}^\phi(e_b,e_t)
+\sum_{b}\eta_{kl}^\phi(e_s,e_b)\eta_{ij}^\phi(e_b,e_t)\\
&=& -\sum_{b}[\hat\eta_{kl}^\phi]_{tb}[\hat\eta_{ij}^\phi]_{bs}
+\sum_{b}[\hat\eta_{ij}^\phi]_{tb}[\hat\eta_{kl}^\phi]_{bs}\\
&=&[\hat\eta_{ij}^\phi,\hat\eta_{kl}^\phi]_{ts}\label{eq:long-from-bracket-calculation}
\end{eqnarray*}
is the entry in row $t$ and column $s$ of the matrix $[\hat\eta_{ij}^\phi,\hat\eta_{kl}^\phi]$.

Secondly, we prove that the endomorphisms $\hat\eta_{kl}^\phi$ satisfy the commutation
relations of $\mathfrak{so}(r)$:
\begin{enumerate}
 \item If $1\leq i,j,k,l\leq r$ are all different,
 \begin{equation}
   [\hat\eta_{kl}^\phi,\hat\eta_{ij}^\phi]=0.\label{eq:[kl,ij]=0}
 \end{equation}
\item If $1\leq i,j,k\leq r$ are all different,
 \begin{equation}
 [\hat\eta_{ij}^\phi,\hat\eta_{jk}^\phi]= -\hat\eta_{ik}^\phi.\label{[ij,jk]=-ik}   
 \end{equation}
\end{enumerate}

To prove \rf{eq:[kl,ij]=0}, note that by \rf{eq:despeje1},
\begin{eqnarray}
\eta_{ij}^\phi\cdot f_kf_l\cdot \phi 
&=&\eta_{kl}^\phi\cdot f_if_j\cdot \phi,\label{eq:identity1}
\end{eqnarray}
by \rf{eq:long-from-bracket-calculation} 
\begin{eqnarray*}
{\rm Re}\left<e_s e_t \cdot \eta_{ij}^\phi\cdot f_kf_l\cdot\phi,\phi\right> 
&=&[\hat\eta_{ij}^\phi,\hat\eta_{kl}^\phi]_{ts},\\
{\rm Re}\left<e_se_t \cdot \eta_{kl}^\phi\cdot f_if_j\cdot\phi,\phi\right> 
&=&[\hat\eta_{kl}^\phi,\hat\eta_{ij}^\phi]_{ts},
\end{eqnarray*}
and by \rf{eq:identity1} and the anticommutativity of the bracket, 
\[[\hat\eta_{ij}^\phi,\hat\eta_{kl}^\phi]=0.\]

To prove \rf{[ij,jk]=-ik}, note that by \rf{eq:despeje1}
\begin{eqnarray*}
f_if_j\cdot\eta_{jk}^\phi\cdot \phi 
&=& f_if_k\cdot\phi
\end{eqnarray*}
and
\begin{eqnarray*}
f_jf_k\cdot\eta_{ij}^\phi\cdot \phi 
&=& -f_if_k\cdot\phi
\end{eqnarray*}
so that
\begin{eqnarray*}
f_jf_k\cdot\eta_{ij}^\phi\cdot \phi 
&=&f_if_j\cdot\eta_{jk}^\phi\cdot \phi  - 2 f_if_k\cdot \phi.
\end{eqnarray*}
Thus, 
\begin{eqnarray*}
{\rm Re}\left<e_s e_t \cdot \eta_{ij}^\phi\cdot f_jf_k\cdot\phi,\phi\right>
&=&  {\rm Re}\left<e_s e_t \cdot \eta_{jk}^\phi\cdot
f_if_j\cdot\phi,\phi\right>- 2 \eta_{ik}^\phi(e_s,e_t)
\end{eqnarray*}
and by \rf{eq:long-from-bracket-calculation}
\begin{eqnarray*}
 [\hat\eta_{ij}^\phi,\hat\eta_{jk}^\phi] &=& [\hat\eta_{jk}^\phi,\hat\eta_{ij}^\phi ] - 2\hat\eta_{ik}^\phi, 
\end{eqnarray*}
i.e.
\begin{eqnarray*}
 [\hat\eta_{ij}^\phi,\hat\eta_{jk}^\phi] &=&  -\hat\eta_{ik}^\phi.
\end{eqnarray*}

Thirdly, we will prove, in five separate cases, that the set of endomorphisms $\{\hat\eta_{kl}^\phi\}$ is
linearly
independent. 
For $r=0,1$ there are no endomorphisms.
For $r=2$ it is obvious since there is only one non-zero endomorphism. 
For $r=3$, suppose
\begin{eqnarray*}
0 &=& \alpha_{12}\hat\eta_{12}^\phi + \alpha_{13}\hat\eta_{13}^\phi+\alpha_{23}\hat\eta_{23}^\phi,
\end{eqnarray*}
where $\alpha_{12}\not=0$.
Take the Lie bracket with $\hat\eta_{13}^\phi$ to get
\begin{eqnarray*}
 0
   &=&  \alpha_{12}\hat\eta_{23}^\phi 
-\alpha_{23}\hat\eta_{12}^\phi ,
\end{eqnarray*}
i.e.
\[\hat\eta_{23}^\phi ={\alpha_{23}\over \alpha_{12}}\hat\eta_{12}^\phi.\]
We can also consider the bracket with $\hat\eta_{23}^\phi$,
\begin{eqnarray*}
 0
   &=&  -\alpha_{12}\hat\eta_{13}^\phi 
+\alpha_{13}\hat\eta_{12}^\phi ,
\end{eqnarray*}
so that
\[\hat\eta_{13}^\phi ={\alpha_{13}\over \alpha_{12}}\hat\eta_{12}^\phi.\]
By substituting in the original equation we get
\begin{eqnarray*}
0 
   &=& (\alpha_{12}^2 + \alpha_{13}^2+\alpha_{23}^2)\hat\eta_{12}^\phi,
\end{eqnarray*}
which gives a contradiction.

Now suppose $r\geq 5$ and that there is a linear combination
\[0=\sum_{k<l} \alpha_{kl}\hat\eta_{kl}^\phi .\]
Taking succesive brackets with $\hat\eta_{13}^\phi$, $\hat\eta_{12}^\phi$, $\hat\eta_{34}^\phi$ and
$\hat\eta_{45}^\phi$ we get the identity
\[\alpha_{12}\hat\eta_{15}^\phi=0,\] 
i.e. $\alpha_{12}=0$. Similar arguments give the vanishing of every $\alpha_{kl}$.

For $r=4$, suppose there is a linear combination
\[0= 
\alpha_{12}\eta_{12}^\phi +
\alpha_{13}\eta_{13}^\phi +
\alpha_{14}\eta_{14}^\phi +
\alpha_{23}\eta_{23}^\phi +
\alpha_{24}\eta_{24}^\phi +
\alpha_{34}\eta_{34}^\phi.
\]
Multiply by $-\phi$
\begin{eqnarray*}
 0
   &=& 
(\alpha_{12}f_{12} +
\alpha_{13}f_{13} +
\alpha_{14}f_{14} +
\alpha_{23}f_{23} +
\alpha_{24}f_{24} +
\alpha_{34}f_{34})\cdot\phi. 
\end{eqnarray*}
Multiply by $-f_{12}$
\begin{eqnarray*}
 0
   &=& 
(\alpha_{12} 
-\alpha_{13}f_{23} 
-\alpha_{14}f_{24} 
+\alpha_{23}f_{13} 
+\alpha_{24}f_{14} 
-\alpha_{34}f_{1234})\cdot\phi. 
\end{eqnarray*}
Now, take hermitian product with $\phi$
\begin{eqnarray*}
 0
   &=& 
\left<(\alpha_{12} 
-\alpha_{13}f_{23} 
-\alpha_{14}f_{24} 
+\alpha_{23}f_{13} 
+\alpha_{24}f_{14} 
-\alpha_{34}f_{1234})\cdot\phi,\phi\right>\\ 
   &=& 
\alpha_{12}|\phi|^2 
-\alpha_{34}\left<f_{1234}\cdot\phi,\phi\right>\\ 
   &=& 
\alpha_{12}.
\end{eqnarray*}
Similar arguments give the vanishing of the other coefficients.
\qd

%

\begin{lemma}\label{lemma:kernel}
Let $r\geq2$ and $\phi\in\Sigma_r\otimes\Delta_n$ be a partially pure spinor. Then
 \[V^\phi\subseteq \bigcap_{1\leq k<l\leq r}\ker\hat\eta_{kl}^\phi.\]
\end{lemma}
{\em Proof}. Let $1\leq k<l \leq r$ be fixed 
and $X\in V^\phi$. 
Since $\mathbb{R}^n=V^\phi \oplus (V^\phi)^\perp$
and $J^\phi$ is a complex structure on $V^\phi$, there exists a basis 
$\{e_1,e_2,\ldots,e_{2m-1},e_{2m}\}\cup\{e_{2m+1},\ldots,e_{2m+r}\}$ 
such that 
\begin{eqnarray*}
 V^\phi &=& {\rm span}(e_1,e_2,\ldots,e_{2m-1},e_{2m}),\\
 (V^\phi)^\perp &=& {\rm span}(e_{2m+1},\ldots,e_{2m+r}),\\
 J^\phi(e_{2j-1}) &=& e_{2j},\\
 J^\phi(e_{2j}) &=& -e_{2j-1},
\end{eqnarray*}
where $m=(n-r)/2$ and $1\leq j \leq m$.
Note that
\begin{eqnarray*}
 \hat\eta_{kl}^\phi(e_{2j-1}) 
   &=& \sum_{a=1}^n{\rm Re}\left< e_{2j-1}\wedge e_a\cdot f_{kl}\cdot\phi,\phi\right>e_a\\
   &=& -\sum_{a\not= 2j-1}^n{\rm Re}\left<  f_{kl}\cdot e_a e_{2j-1} \cdot\phi,\phi\right>e_a\\
   &=& -\sum_{a\not= 2j-1}^n{\rm Re}\left<  f_{kl}\cdot e_a (iJ^\phi (e_{2j-1})) \cdot\phi,\phi\right>e_a\\
   &=& \sum_{a\not= 2j-1}^n{\rm Im} \left<   e_a e_{2j} \cdot f_{kl}\cdot\phi,\phi\right>e_a\\
   &=& -{\rm Im}\left<f_{kl}\cdot\phi,\phi\right>e_{2j}\\
   &=&0.
\end{eqnarray*}
\qd

\begin{lemma}
Let $r\geq2$ and $\phi\in\Sigma_r\otimes\Delta_n$ be a partially pure spinor.
Then $(V^\phi)^\perp$ carries a standard representation of $\mathfrak{so}(r)$, and an orientation.
\end{lemma}
{\em Proof}.
By Lemma \ref{lemma:lie-algebra}, $\mathfrak{so}(r)$ is represented non-trivially on $\mathbb{R}^n =
V^\phi \oplus (V^\phi)^\perp$ and, by Lemma \ref{lemma:kernel}, it acts trivially on $V^\phi$.
Thus $(V^\phi)^\perp$ is a nontrivial representation of $\mathfrak{so}(r)$ of dimension $r$.
\qd

{\bf Remark}.
The existence of a partially pure spinor implies $r\equiv n$ (mod $2$).
In this case, let
$(e_1,\ldots, e_n)$ and $(f_1,\ldots, f_r)$ be orthonormal frames for $\mathbb{R}^n$ and
$\mathbb{R}^r$ respectively, 
\[{\rm vol}_n= e_1\cdots e_n, \quad\quad {\rm vol}_r= f_1\cdots f_r,\]
and 
\begin{eqnarray*}
F:\Sigma_r\otimes\Delta_n &\lra& \Sigma_r\otimes\Delta_n\\ 
\phi &\mapsto& (-i)^{n/2}i^{r/2}{\rm vol}_n\cdot {\rm vol}_r  \cdot\phi.
\end{eqnarray*}
Note that $i^{r/2}{\rm vol}_r$ acts as $(-1)^{r/2} {\rm Id}_{\Sigma_r}$ on $\Sigma_r$ and that
$(-i)^{n/2}{\rm vol}_n$ determines the decomposition $\Delta_n=\Delta_n^+\oplus \Delta_n^-$.
Thus we have that
\[\Sigma_r\otimes\Delta_n = (\Sigma_r\otimes\Delta_n)^+ \oplus (\Sigma_r\otimes\Delta_n)^-,\]
and we will call elements in $(\Sigma_r\otimes\Delta_n)^+$ and $(\Sigma_r\otimes\Delta_n)^-$
positive and negative twisted spinors respectively.

\begin{defi} Let $n$ be even, $\mathbb{R}^n$ be endowed with the standard inner product and orientation, and
${\rm vol}_n$ denote the volume form. Let
$V$, $W$ be two orthogonal oriented subspaces such that $\mathbb{R}^n = V\oplus W$. Furthermore, assume $V$
admits a complex structure inducing the given orientation on $V$. The oriented triple $(V,J,W)$ will be called
{\em positive} if given (oriented) orthonormal frames $(v_1,J(v_1),\ldots,v_m,J(v_m))$ and $(w_1,\ldots,w_r)$
of $V$ and $W$ respectively, 
\[v_1\wedge J(v_1)\wedge \ldots\wedge v_m\wedge J(v_m)\wedge w_1\wedge\ldots\wedge w_r = {\rm vol}_n,\]
and {\em negative} if
\[v_1\wedge J(v_1)\wedge \ldots\wedge v_m\wedge J(v_m)\wedge w_1\wedge\ldots\wedge w_r = -{\rm vol}_n.\]
\end{defi}

\begin{lemma}
 If $r$ is even, a partially pure spinor $\phi$ is either positive or negative. Furthermore, 
a partialy pure spinor $\phi$ is positive (resp. negative) if and only if 
the corresponding oriented triple $(V^\phi, J^\phi, (V^\phi)^\perp)$ is positive (resp. negative).
\end{lemma}
{\em Proof}.
We must prove that either
$\phi\in(\Sigma_r\otimes\Delta_n)^+$ or
$\phi\in(\Sigma_r\otimes\Delta_n)^-$.
Since $\phi$ is a partially pure spinor, 
there exist frames
$(e_1',\ldots, e_{2m}')$ and $(e_{2m+1}',\ldots,e_{2m+r}')$ 
of $V^\phi$ and $(V^\phi)^\perp$ respectively such that
\[e_{2j}'= J(e_{2j-1}')\quad\mbox{and}\quad 
\eta_{kl}^\phi = e_{2m+k}'\wedge e_{2m+l}',\]
where $1\leq j\leq m$ and $1\leq k<l\leq r$.
Now, 
\[e_1'\wedge e_2'\wedge \ldots \wedge e_{2m}'\wedge e_{2m+1}'\wedge\ldots\wedge e_{2m+r}'=\pm {\rm vol}_n.\]
Then,
\begin{eqnarray*}
 (-i)^{n/2}i^{r/2}{\rm vol}_n\cdot {\rm vol}_r  \cdot\phi 
   &=& 
\pm  (-i)^{n/2}i^{r/2} e_1' e_2' \cdots  e_{2m}' e_{2m+1}'\cdots e_{2m+r}'\cdot f_1\cdots f_r\cdot \phi\\   
   &=& 
\pm  (-i)^{n/2}i^{r/2} e_1' J(e_1') \cdots  e_{2m-1}'J(e_{2m-1}') \eta_{12}^\phi\cdots
\eta_{r-3,r-2}^\phi\cdot
f_{12}\cdots f_{r-1,r}\cdot \eta_{r-1,r}^\phi\cdot\phi\\  
   &=& 
\pm  (-i)^{n/2}i^{r/2} e_1' J(e_1') \cdots  e_{2m-1}'J(e_{2m-1}') \eta_{12}^\phi\cdots
\eta_{r-3,r-2}^\phi\cdot
f_{12}\cdots f_{r-3,r-2}\cdot\phi\\  
   &=& 
\pm  (-i)^{n/2}i^{r/2} e_1' J(e_1') \cdots  e_{2m-1}'J(e_{2m-1}') \cdot\phi\\  
   &=& 
\pm  (-i)^{n/2}i^{r/2} e_1' J(e_1') \cdots  e_{2m-3}'J(e_{2m-3}')e_{2m-1}'(-ie_{2m-1}') \cdot\phi\\  
   &=& 
\pm   (-1)^m(-i)^{n/2+m}i^{r/2} \phi\\  
   &=& 
 \pm   \phi,  
\end{eqnarray*}
i.e.
$\phi\in(\Sigma_r\otimes\Delta_n)^\pm$.
\qd

\subsection{Orbit of a partially pure spinor}

\begin{lemma} \label{lemma:g(phi)-partially-pure-spinor}
Let $\phi\in\Sigma_r\otimes\Delta_n$ be a
partially pure spinor. If $g\in Spin^{c,r}(n)$, then $g(\phi)$ is also a partially pure spinor.
\end{lemma}
{\em Proof}.
Let $g\in Spin^{c,r}(n)$ and 
$\lambda_n^{c,r}(g)=(g_1,g_2,g_3)\in SO(n)\times SO(r)\times U(1)$.
First, suppose $X,Y\in V^\phi$, 
\[X\cdot \phi = i \,\, Y\cdot \phi.\]
Apply $g$ on both sides 
\begin{eqnarray*}
g_1(X) \cdot g(\phi)
   &=& i\,\,g_1(Y)\cdot g(\phi). 
\end{eqnarray*}
which means that $g_1$ maps $V^\phi$ into $V^{g(\phi)}$ injectively. 
On the other hand, any pair of
vectors $\tilde X, \tilde Y\in V^{g(\phi)}$ such that
\[\tilde X\cdot g(\phi) = i \,\, \tilde Y\cdot g(\phi),\]
are the image under $g_1$ of some vectors $X,Y\in\mathbb{R}^n$, i.e.
\[g_1(X)\cdot g(\phi) = i \,\, g_1(Y)\cdot g(\phi).\]
Apply $g^{-1}$ on both sides to get
\[X\cdot \phi = i\,\, Y\cdot \phi,\]
so that $X,Y\in V^\phi$, i.e. 
$V^{g(\phi)} = g_1(V^\phi)$.
Moreover, 
\[ J^{g(\phi)} = g_1|_{V^\phi}\circ J^\phi \circ (g_1|_{v^\phi})^{-1}.\]
Now, let $e_a'=g_1^{-1}(e_a)$  and $f_k'=g_2^{-1}(f_k)$, so that
\begin{eqnarray*}
 \eta_{kl}^{g(\phi)}\cdot g(\phi)
   &=& \sum_{1\leq a<b\leq n}  \eta_{kl}^{g(\phi)}(e_a,e_b)e_ae_b\cdot g(\phi) \\
   &=& \sum_{1\leq a<b\leq n}  \left<g_1(e_a')g_1(e_b') \cdot g_2(f_k')g_2(f_l')\cdot g(\phi),
g(\phi)\right>g_1(e_a')g_1(e_b')\cdot g(\phi) \\
   &=& \sum_{1\leq a<b\leq n} \left<e_a'e_b' \cdot f_k'f_l'\cdot \phi, \phi\right>g(e_a'e_b'\cdot \phi)\\
   &=& g\left(\sum_{1\leq a<b\leq n} \eta_{kl}'^{\phi}(e_a',e_b')e_a'e_b'\cdot\phi\right)\\
   &=& g\left(\eta_{kl}'^{\phi}\cdot\phi\right)\\
   &=& g\left(-f_k'f_l'\cdot\phi\right)\\
   &=& -f_kf_l\cdot g(\phi),
\end{eqnarray*}
and
\begin{eqnarray*}
 \left< f_kf_l\cdot g(\phi),g(\phi)\right>
   &=& \left< g(f_k'f_l'\cdot \phi), g(\phi)\right>\\
   &=& \left< f_k'f_l'\cdot \phi, \phi\right>\\
   &=& 0.
\end{eqnarray*}
For $r=4$, note that the volume form is invariant under
$SO(4)$ 
\[\left<f_1f_2f_3f_4\cdot g(\phi),g(\phi)\right> = \left<f_1f_2f_3f_4\cdot \phi,\phi\right> = 0.\]
\qd

\begin{lemma}
Let $\phi\in\Sigma_r\otimes\Delta_n$ be a partially pure spinor. 
The stabilizer of $\phi$ is isomorphic to $U(m)\times SO(r)$.
\end{lemma}
{\em Proof}. Let $g\in Spin^{c,r}(n)$ be such that
$g(\phi)=\phi$ and $\lambda_n^{c,r}(g)=(g_1,g_2,g_3)\in SO(n)\times SO(r)\times U(1)$.
It can be checked easily that
\begin{eqnarray*}
 [g_1,J^\phi] &=& 0 \\
 g_1(V^\phi) &=& V^\phi \\
 g_1|_{V^\phi} &\in& U(V^\phi,J^\phi) \cong U(m).
\end{eqnarray*}
%
Clearly, $g_1((V^\phi)^\perp)= (V^\phi)^\perp$.

As in Lemma \ref{lemma:kernel}, one can prove
\[\eta_{kl}^\phi = \sum_{2m+1\leq a< b \leq 2m+r} \eta_{kl}^\phi(e_a,e_a)e_ae_b\in \ext^2(V^\phi)^\perp,\]
where $(e_1,\ldots,e_{2m+r})$ is an oriented frame of $V^\phi \oplus (V^\phi)^\perp$.
Furthermore,
\begin{eqnarray*}
 g_1(\eta_{kl}^\phi)
   &=& 
   \eta_{kl}'^\phi,
\end{eqnarray*}
where 
$f_k'=g_2(f_k)$.
Now, we have that
\[
\begin{array}{ccc}
f_kf_l & \xrightarrow{g_2} & f_k'f_l'\\
\downarrow &  & \downarrow\\
\eta_{kl}^\phi & \xrightarrow{h_2} & \eta_{kl}'^\phi
\end{array}
\]
for the diagram
\[
\begin{array}{ccccccc}
\mathfrak{so}(r)&\cong&\ext^2\mathbb{R}^r & \xrightarrow{g_2} & \ext^2\mathbb{R}^r&\cong&\mathfrak{so}(r)\\
&&\downarrow &  & \downarrow&&\\
\mathfrak{so}(r)&\cong&\ext^2(V^\phi)^\perp & \xrightarrow{h2} & \ext^2(V^\phi)^\perp &\cong& \mathfrak{so}(r)
\end{array}
\]
where the vertical arrows are Lie algebra isomorphisms and the horizontal arrows correspond to $g_2$ and $h_2$
acting
via the adjoint representation of $SO(r)$. Thus, $h_2$ and $g_2$ correspond to each other under the
isomorphism 
$\ext^2(V^\phi)^\perp\cong \ext^2\mathbb{R}^r$ given by
$f_{kl}\mapsto \eta_{kl}^\phi$.

Since $h_1$ is unitary with respect to $J$, there is a frame $(e_1,\ldots, e_{2m})$ of $V^\phi$ such that
\[e_{2j}= J(e_{2j-1})\]
and $h_1$ is diagonal with respect to the unitary basis
$\{e_{2j-1}-ie_{2j}\,|\, j=1,\ldots,m\}$,
i.e.
\[h_1(e_{2j-1}-ie_{2j})=e^{i\theta_j}(e_{2j-1}-ie_{2j})\]
where $0\leq \theta_j<2\pi$.
On the other hand, 
there is a frame $(f_1,\ldots, f_r)$ of $\mathbb{R}^r$ such that
\[g_2= R_{\varphi_1}\circ \cdots \circ R_{\varphi_{[r/2]}}
\]
where $R_{\varphi_k}$ is a rotation by an angle $\varphi_k$ on the plane generated by 
$f_{2k-1}$ and $f_{2k}$, $1\leq k\leq [r/2]$. Now, since the endomorphisms $\hat\eta_{kl}^\phi$ span an
isomorphic
copy of $\mathfrak{so}(r)$, there is a frame 
$(e_{2m+1},\ldots,e_{2m+r})$ of $(V^\phi)^\perp$ such that
\[\eta_{kl}^\phi = e_{2m+k}\wedge e_{2m+l},\]
$1\leq k<l\leq r$.
Since the adjoint representation of $SO(r)$ is faithful
\[h_2= R_{\varphi_1}'\circ \cdots \circ R_{\varphi_{[r/2]}}'
\]
where $R_{\varphi_k}'$ is a rotation by an angle $\varphi_k$ on the plane generated by 
$e_{2m+2k-1}$ and $e_{2m+2k}$, $1\leq k\leq [r/2]$.
Thus, 
\begin{eqnarray*}
g
   &=&
 \pm\left[\prod_{j=1}^m (\cos(\theta_j/2)- \sin(\theta_j/2)e_{2j-1}e_{2j})\cdot
\prod_{k=1}^{[r/2]} (\cos(\varphi_k/2)- \sin(\varphi_k/2)\eta_{2k-1,2k}^\phi),\right.\\
   &&
 \left. 
\prod_{k=1}^{[r/2]} (\cos(\varphi_k/2)- \sin(\varphi_k/2)f_{2k-1}f_{2k})
, e^{i\theta/2}
\right] .
\end{eqnarray*}
Now, 
\begin{eqnarray*}
 \phi 
   &=&
  g(\phi)\\
   &=&
 \pm e^{i\theta/2}\prod_{j=1}^m (\cos(\theta_j/2)- \sin(\theta_j/2)e_{2j-1}e_{2j})
\\
   &&
\cdot\prod_{k=1}^{[r/2]} (\cos(\varphi_k/2)- \sin(\varphi_k/2)\eta_{2k-1,2k}^\phi)\cdot\prod_{k=1}^{[r/2]}
(\cos(\varphi_k/2)- \sin(\varphi_k/2)f_{2k-1}f_{2k})  (\phi)\\
   &=&
 \pm e^{i\theta/2}\prod_{j=1}^m (\cos(\theta_j/2)- \sin(\theta_j/2)e_{2j-1}e_{2j})
\cdot\prod_{k=1}^{[r/2]-1} (\cos(\varphi_k/2)- \sin(\varphi_k/2)\eta_{2k-1,2k}^\phi)\\
   &&
\cdot\prod_{k=1}^{[r/2]}
(\cos(\varphi_k/2)- \sin(\varphi_k/2)f_{2k-1}f_{2k})
(\cos(\varphi_{[r/2]}/2)-
\sin(\varphi_{[r/2]}/2)\eta_{2{[r/2]}-1,2{[r/2]}}^\phi)\cdot  (\phi)\\
   &=&
 \pm e^{i\theta/2}\prod_{j=1}^m (\cos(\theta_j/2)- \sin(\theta_j/2)e_{2j-1}e_{2j})
\\
   &&
\cdot\prod_{k=1}^{[r/2]-1} (\cos(\varphi_k/2)- \sin(\varphi_k/2)\eta_{2k-1,2k}^\phi)\cdot
\prod_{k=1}^{[r/2]-1}
(\cos(\varphi_k/2)- \sin(\varphi_k/2)f_{2k-1}f_{2k})\\
   &&
(\cos(\varphi_{[r/2]}/2)- \sin(\varphi_{[r/2]}/2)f_{2{[r/2]}-1}f_{2{[r/2]}})\cdot
(\cos(\varphi_{[r/2]}/2)+\sin(\varphi_{[r/2]}/2)f_{2{[r/2]}-1}f_{2{[r/2]}})\cdot  (\phi)\\
   &=&
 \pm e^{i\theta/2}\prod_{j=1}^m (\cos(\theta_j/2)- \sin(\theta_j/2)e_{2j-1}e_{2j})
\\
   &&
\cdot\prod_{k=1}^{[r/2]-1} (\cos(\varphi_k/2)- \sin(\varphi_k/2)\eta_{2k-1,2k}^\phi)\cdot
\prod_{k=1}^{[r/2]-1}
(\cos(\varphi_k/2)- \sin(\varphi_k/2)f_{2k-1}f_{2k})  (\phi)\\
   &=&
 \pm e^{i\theta/2}\prod_{j=1}^m (\cos(\theta_j/2)- \sin(\theta_j/2)e_{2j-1}e_{2j})
 (\phi)\\
   &=&
 \pm e^{i\theta/2}\prod_{j=1}^m (\cos(\theta_j/2)+i \sin(\theta_j/2)e_{2j-1}(iJ(e_{2j-1})))
 (\phi)\\
   &=&
 \pm e^{i\theta/2}\prod_{j=1}^m (\cos(\theta_j/2)+i \sin(\theta_j/2)e_{2j-1}e_{2j-1})
 (\phi)\\
   &=&
 \pm e^{i\theta/2}\prod_{j=1}^m (\cos(\theta_j/2)-i \sin(\theta_j/2))
 (\phi)\\
   &=&
 \pm e^{i\theta/2}\prod_{j=1}^m e^{-i\theta_j/2}
 (\phi)\\
   &=&
 \pm e^{{i\over 2}(\theta-\sum_{j=1}^m\theta_j)}
 (\phi).
\end{eqnarray*}
This means
\[e^{{i\over 2}(\theta-\sum_{j=1}^m\theta_j)}=\pm1\]
i.e.
\begin{eqnarray*}
{\det}_{\mathbb{C}}(h_1) 
   &=& 
  e^{i\sum_{j=1}^m \theta_j} \\
   &=& 
  e^{i\theta} .
\end{eqnarray*}
Thus we have found that
\[\lambda_{n}^{c,r}(g) =  ((h_1,h_2),h_2,{\det}_{\mathbb{C}}(h_1)) ,\]
which is in the image of the horizontal row in the diagram of Lemma \ref{lemma:subgroup2}
\[
\begin{array}{ccc}
 &  & Spin^{c,r}(n)\\
 & \nearrow & \downarrow\\
U(m)\times SO(r) & \rightarrow  & SO(n)\times SO(r)\times U(1)
\end{array}
\]
\qd

{\bf Remark}.
Note that for any spinor $\phi\in\Sigma_r\otimes\Delta_n$, $g\in Spin^{c,r}(n)$, $\lambda_n^{c,r}(g)\in
SO(n)\times SO(r)\times U(1)$,
\begin{eqnarray*}
\eta_{kl}^{g(\phi)}(X,Y) 
   &=&
 \left<X\wedge Y\cdot f_kf_l \cdot g(\phi),g(\phi)\right> \\
   &=&  \left<g_1(X')\wedge g_1(Y')\cdot g_2(f_k')g_2(f_l') \cdot g(\phi),g(\phi)\right> \\
   &=&  \left<g(X'\wedge Y'\cdot f_k'f_l' \cdot \phi),g(\phi)\right> \\
   &=&  \left<X'\wedge Y'\cdot f_k'f_l' \cdot \phi,\phi\right> \\
   &=:&  \eta_{kl}'^{\phi} (X',Y'),
\end{eqnarray*}
for $X'=g_1^{-1}(X),Y'=g_1^{-1}(Y)\in\mathbb{R}^n$,
$f_k'=g_2^{-1}(f_k)$. Thus, the matrices representing
$\eta_{kl}^{g(\phi)}$ (with respect to some basis) are conjuagte to the matrices representing
$\eta_{kl}'^\phi$.

\begin{lemma}\label{lemma: little orbit}
Let $\phi,\psi\in\Sigma_r\otimes\Delta_n$ be partially pure spinors and
$Spin^c(r)$ the standard copy of this group in $Spin^{c,r}(n)$. 
Then,  $\psi\in Spin^{c}(r)\cdot\phi$ if and only if they generate the same oriented tiple
$(V^\phi,J^\phi,(V^\phi)^\perp)=(V^\psi,J^\psi,(V^\psi)^\perp)$.
\end{lemma}
{\em Proof}. Suppose $\psi = g(\phi)$ for some $g\in Spin^c(r)\subset Spin^{c,r}(n)$, and
let $\lambda_n^{c,r}(g)=(1,g_2,e^{i\theta})$.
Such an element induces 
\begin{eqnarray*}
 \left<X\wedge Y\cdot f_kf_l \cdot g(\phi),g(\phi)\right>
   &=&  \left<X\wedge Y\cdot f_k'f_l' \cdot \phi,\phi\right> 
\end{eqnarray*}
for $f_k'=g_2^{-1}(f_k)$, i.e.
\[\eta_{kl}^{g(\phi)}(X,Y) =  \eta_{kl}'^{\phi}(X,Y),
\]
so that they span the same copy of $\mathfrak{so}(r)$ in $\End^-(\mathbb{R}^n)$, 
\[ {\rm span}(\eta_{kl}^{g(\phi)}) = {\rm span}(\eta_{kl}'^\phi)\cong
\mathfrak{so}(r)\subset\End^-(\mathbb{R}^n).
\]
Thus, by Lemma \ref{lemma:g(phi)-partially-pure-spinor},
the partially pure spinors $\phi$ and $g(\phi)$ determine the same oriented triple 
$(V^{g(\phi)}, J^{g(\phi)},(V^{g(\phi)})^\perp ) =(V^\phi,V^\phi, (V^\phi)^\perp)$.

Conversely, assume $(V^\phi,J^\phi,(V^\phi)^\perp)=(V^\psi,J^\psi,(V^\psi)^\perp)$, and
consider the subalgebras of
\begin{eqnarray*}
\mathfrak{so}(r)^\phi &=& {\rm span}(\eta_{kl}^\phi + f_{kl})\\
\mathfrak{so}(r)^\psi &=& {\rm span}(\eta_{kl}^\psi + f_{kl}).
\end{eqnarray*}
There exist frames $(e_{2m+1},\ldots,e_{2m+r})$ and $(e_{2m+1}',\ldots,e_{2m+r}')$ of
$(V^\phi)^\perp$ and $ (V^\psi)^\perp$ respectively, such that
\begin{eqnarray*}
 \eta_{kl}^\phi &=& e_{2m+k}\wedge e_{2m+l},\\
 \eta_{kl}^\psi &=& e_{2m+k}'\wedge e_{2m+l}'.
\end{eqnarray*}
Let $A=(a_{kl})\in SO(r)$ the matrix such that 
\[A:\,\,e_{2m+k}'\mapsto a_{k1}e_{2m+1}' + \cdots + a_{kr}e_{2m+r}' = e_{2m+k}\]
$1\leq k < l \leq r$.
The induced transformation maps
\[A:\,\,e_{2m+k}'\wedge e_{2m+l}' \mapsto e_{2m+k}\wedge e_{2m+l},\]
and set
\[A^T:\,\,  f_k\mapsto a_{1k}f_{1} + \cdots + a_{rk}f_{r}=f_k',\]
and 
\[A^T:\,\,f_k\wedge f_l\mapsto f_k'\wedge f_l'.\]
Consider
\begin{eqnarray*}
\left<e_{2m+p}\wedge e_{2m+q}\cdot f_k'f_l'\cdot \psi,\psi\right> 
   &=&
  \left<\left(\sum_{s=1}^r a_{ps}e_{2m+s}'\right)\wedge \left(\sum_{t=1}^r a_{qt}e_{2m+t}'\right)\cdot
f_k'f_l'\cdot \psi,\psi\right>\\ 
   &=&
  \left<\left(\sum_{s<t} (a_{ps}a_{qt}-a_{pt}a_{qs})e_{2m+s}'\wedge e_{2m+t}'\right)\cdot
f_k'f_l'\cdot \psi,\psi\right>\\ 
   &=&
  \sum_{s<t} (a_{ps}a_{qt}-a_{pt}a_{qs})\left<e_{2m+s}'\wedge e_{2m+t}'\cdot
f_k'f_l'\cdot \psi,\psi\right>\\ 
   &=&
  \sum_{s<t} (a_{ps}a_{qt}-a_{pt}a_{qs})\left<e_{2m+s}'\wedge e_{2m+t}'\cdot
\left(\sum_{i=1}^r a_{ik}f_i\right)\left(\sum_{j=1}^r a_{jl}f_j\right)\cdot \psi,\psi\right>\\ 
   &=&
  \sum_{s<t} (a_{ps}a_{qt}-a_{pt}a_{qs})\left<e_{2m+s}'\wedge e_{2m+t}'\cdot
\left(\sum_{i<j} (a_{ik}a_{jl}-a_{il}a_{jk})f_if_j\right)\cdot \psi,\psi\right>\\ 
   &=&
  \sum_{s<t}\sum_{i<j} (a_{ps}a_{qt}-a_{pt}a_{qs})(a_{ik}a_{jl}-a_{il}a_{jk})\left<e_{2m+s}'\wedge
e_{2m+t}'\cdot
 f_if_j\cdot \psi,\psi\right>\\ 
   &=&
  \sum_{s<t}\sum_{i<j} (a_{ps}a_{qt}-a_{pt}a_{qs})(a_{ik}a_{jl}-a_{il}a_{jk}) \delta_{si}\delta_{tj}\\
   &=&
  \sum_{s<t} (a_{ps}a_{qt}-a_{pt}a_{qs})(a_{sk}a_{tl}-a_{sl}a_{tk}) \\
   &=&
  \delta_{pk}\delta_{ql} ,
\end{eqnarray*}
since the induced tranformation by $A$ on $\ext^2\mathbb{R}^r$ is orthogonal.
This means
\[\eta_{kl}'^\psi = \eta_{kl}^\phi= e_{2m+k}\wedge e_{2m+l}.\]
Now consider $g\in Spin^c(r)\subset Spin^{c,r}(n)$ such that $\lambda_n^{c,r}(g)=(1,A,1)\in SO(n)\times
SO(r)\times U(1)$.
Then
\begin{eqnarray*}
\delta_{pk}\delta_{ql}
   &=& 
\left<e_{2m+p}\wedge e_{2m+q}\cdot f_k'f_l'\cdot \psi,\psi\right> \\
   &=& 
\left<g(e_{2m+p}\wedge e_{2m+q}\cdot f_k'f_l'\cdot \psi),g(\psi)\right> \\
   &=& 
\left<e_{2m+p}\wedge e_{2m+q}\cdot A(f_k')A(f_l')\cdot g(\psi),g(\psi)\right> \\
   &=& 
\left<e_{2m+p}\wedge e_{2m+q}\cdot f_kf_l\cdot g(\psi),g(\psi)\right> ,
\end{eqnarray*}
i.e.
\[\eta_{kl}^{g(\psi)} = e_{2m+k}\wedge e_{2m+l} = \eta_{kl}^\phi,\]
so that
\[\mathfrak{so}(r)^{g(\psi)} = {\rm span}(\eta_{kl}^{g(\psi)} + f_{kl}) = {\rm span}(\eta_{kl}^{\phi} +
f_{kl})=\mathfrak{so}(r)^\phi.\]
This implies that $g(\psi)$ and $\phi$ share the same stabilizer
\[  U(V^\phi,J^\phi)\times exp(\mathfrak{so}(r)^\phi)= U(V^{g(\psi)},J^{g(\psi)})\times
exp(\mathfrak{so}(r)^{g(\psi)}) \cong U(m)\times SO(r).\]
But there is only a 1-dimensional summand in the decomposition of $\Sigma_r\otimes \Delta_n$
under this subgroup. More precisely, under this subgroup
\begin{eqnarray*}
 \Sigma_r\otimes \Delta_n
   &=&
  \Sigma_r\otimes\Delta_r\otimes\Delta_{2m},
\end{eqnarray*}
where $\Delta_{2m}$ decomposes under $U(m)$ and contains only a $1$-dimensional trivial summand
\cite{Friedrich}, while 
$\Sigma_r\otimes\Delta_r$ is isomorphic to a subspace of the complexified space of alternating forms on
$\mathbb{R}^r$ which also contains only a $1$-dimensional trivial summand.
Thus, $g(\psi)= e^{i\theta} \phi$ for some $\theta\in[0,2\pi)\subset\mathbb{R}$.

\qd

\begin{lemma} 
\begin{itemize}
 \item  If $r$ is odd, the group $Spin^{c,r}(n)$ acts transtitively on the set of partially pure spinors in
$\Sigma_r\otimes\Delta_n $.
 \item  If $r$ is even, the group $Spin^{c,r}(n)$ acts transtitively on the set of positive partially pure
spinors in
$(\Sigma_r\otimes\Delta_n)^+ $.
\end{itemize}
\end{lemma}
{\em Proof}. 
Suppose that $r$ is odd.
Note that the standard partially pure spinor $\phi_0$ satisfies the conditions
\begin{equation}
 \left\{
\begin{array}{rcl}
e_{2j-1}e_{2j}\cdot\phi_0 & = & i\phi_0, \\
e_{2m+k}e_{2m+l}\cdot\phi & = & -f_{kl}\cdot\phi, \\
\left<f_{kl}\cdot\phi,\phi\right> & = & 0, 
\end{array}
\right.\label{eq: equations standard partially pure spinor}
\end{equation}
where $(e_1,\ldots, e_{n})$ and $(f_1,\ldots,f_r)$ are the standard oriented frames of $\mathbb{R}^n$ and
$\mathbb{R}^r$ respectively. 
There exist frames
$(e_1',\ldots, e_{2m}')$ and $(e_{2m+1}',\ldots,e_{2m+r}')$ of $V^\phi$ and $(V^\phi)^\perp$ respectively such
that
\[e_{2j}'= J(e_{2j-1}')\quad
\mbox{and}
\quad\eta_{kl}^\phi = e_{2m+k}'\wedge e_{2m+l}',\]
$1\leq k<l\leq r$, $1\leq j\leq m$.
Call $g_1'\in O(n)$
the transformation of $\mathbb{R}^n$ taking the new frame to the standard one.
Define $g_1\in SO(n)$ as follows
\[\left\{
 \begin{array}{ll}
g_1 = g_1', & \mbox{if $e_1'\wedge\ldots\wedge e_{2m+r}'= {\rm vol}_n$,} \\
g_1 = -g_1', & \mbox{if $e_1'\wedge\ldots\wedge e_{2m+r}'= -{\rm vol}_n$.} \\
 \end{array}
\right.
\]
Then $(g_1,1,1)\in SO(n)\times SO(r) \times U(1)$ has two preimages $\pm\tilde g\in Spin^{c,r}(n)$.
By Lemma \ref{lemma:g(phi)-partially-pure-spinor}, $\tilde{g}(\phi)$ is a partially pure spinor.
We will check that $\tilde{g}(\phi)$ satisfies  \rf{eq: equations standard partially pure
spinor} as $\phi_0$ does.
Indeed,
\begin{eqnarray*}
 e_{2j-1}e_{2j}\cdot \tilde{g}(\phi)
   &=& 
  g_1'(e_{2j-1}')g_1'(e_{2j}')\cdot \tilde{g}(\phi)\\  
   &=& 
  (\pm g_1(e_{2j-1}'))(\pm g_1(e_{2j}'))\cdot \tilde{g}(\phi)\\  
   &=& 
  g_1(e_{2j-1}') g_1(e_{2j}')\cdot \tilde{g}(\phi)\\  
   &=& 
  \tilde{g}(e_{2j-1}'e_{2j}'\cdot \phi)\\  
   &=& 
  \tilde{g}(i \phi)\\  
   &=& 
  i\tilde{g}( \phi), 
\end{eqnarray*}
and
\begin{eqnarray*}
 e_{2m+k}e_{2m+l}\cdot \tilde{g}(\phi)
   &=&
 g_1'(e_{2m+k}')g_1'(e_{2m+l}')\cdot \tilde{g}(\phi)\\
   &=&
 (\pm g_1(e_{2m+k}'))(\pm g_1(e_{2m+l}'))\cdot \tilde{g}(\phi)\\
   &=&
 g_1(e_{2m+k}')g_1(e_{2m+l}')\cdot \tilde{g}(\phi)\\
   &=&
 \tilde{g}(e_{2m+k}'e_{2m+l}'\cdot \phi)\\
   &=&
 \tilde{g}(-f_{k}f_{l}\cdot \phi)\\
   &=&
 -\lambda_2(\tilde{g})(f_{k})\lambda_2(\tilde{g})(f_{l})\cdot \tilde{g}(\phi)\\
   &=&
 -f_{k}f_{l}\cdot \tilde{g}(\phi),
\end{eqnarray*}
since $\lambda_2(\tilde{g})=1$.
Similarly,
\begin{eqnarray*}
 \left< f_kf_l\cdot \tilde{g}(\phi),\tilde{g}(\phi)\right>
   &=&
  \left< \lambda_2(\tilde{g})(f_k)\lambda_2(\tilde{g})(f_l)\cdot \tilde{g}(\phi),\tilde{g}(\phi)\right> \\ 
   &=&
  \left< \tilde{g}(f_kf_l\cdot\phi),\tilde{g}(\phi)\right> \\
   &=&
  \left< f_kf_l\cdot\phi,\phi\right> \\
   &=& 
   0.
\end{eqnarray*}
Thus, $\tilde{g}(\phi)$ generates the same oriented triple $(V^{\tilde g(\phi)},J^{\tilde g(\phi)},
(V^{\tilde g(\phi)})^\perp)=(V^{\phi_0},J^{\phi_0},
(V^{\phi_0})^\perp)$ as $\phi_0$ which, by Lemma \ref{lemma: little orbit}, concludes the proof for $r$ odd. 

The case for $r$ even is similar.
\qd

\begin{theo} \label{theo:characterization}
Let $\mathbb{R}^n$ be endowed with the standard inner product and orientation.
Given $r\in\mathbb{N}$ such that $r< n$, the following objects are equivalent:
\begin{enumerate}
 \item A (positive) triple consisting of a codimension $r$ vector subspace endowed with an orthogonal complex
structure
and an oriented orthogonal complement.
 \item An orbit $Spin^c(r)\cdot\phi$ for some (positive) twisted partially pure spinor $\phi\in
\Delta_n\otimes
\Sigma_r$.
\end{enumerate}
\end{theo}
{\em Proof}. Given a codimension $r$ vector subspace $D$ endowed with an orthogonal complex structure,
$\dim_{\mathbb{R}}(D)=2m$, $n=2m+r$. By Lemma \ref{factorization}
\[\Delta_n \cong \Delta(D^\perp)\otimes\Delta(D) .\]
Let us define
\begin{eqnarray*}
\Sigma_r&\cong& \left\{ 
\begin{array}{ll}
\Delta(D^\perp) & \mbox{if $r$ is odd,}\\
\Delta(D^\perp)^+ & \mbox{if $r$ is even,} 
\end{array}
\right. 
\end{eqnarray*}
so that
\[\Sigma_r\otimes\Delta_n \]
contains the standard twisted partially pure spinor $\phi_0$ of Lemma \ref{lemma:existence}.

The proof of the converse is the content of Subsection   \ref{sec: basic properties}.
\qd

Let $\tilde{\mathcal{S}}$ denote the set of all partially pure spinors of rank $r$
\[\tilde{\mathcal{S}} = {Spin^{c,r}(n) \over U(m)\times SO(r)}.\] 
Consider
\[ \mathcal{S}={\tilde{\mathcal{S}}\over Spin^c(r) } \]
where $Spin^c(r)$ is the canonical copy of such a group in $Spin^{c,r}(n)$.
Thus we have the following expected result.

\begin{corol}
 The space parametrizing subspaces with orthogonal complex structures of codimension $r$ in
$\mathbb{R}^n$ with oriented orthogonal complements is
\[\mathcal{S}\cong{SO(n)\over U(m)\times SO(r)}.\]
\end{corol}
\qd

{\small
\renewcommand{\baselinestretch}{0.5}
\newcommand{\bi}{\vspace{-.05in}\bibitem} }

\end{document}